\documentclass[10pt]{article}
\usepackage[intlimits]{amsmath}
\usepackage{amssymb}
\usepackage{amsmath}
\usepackage{bbm}
\usepackage{epsfig}
\usepackage{float}
\usepackage{cite}
\newtheorem{theorem}{Theorem}[section]

\newtheorem{definition}{Definition}[section]

\newtheorem{proposition}{Proposition}[section]
\newtheorem{lemma}{Lemma}[section]

\newtheorem{rule-def}[theorem]{Rule}
\newtheorem{example}{Example}[section]

\begin{document}
\title{\textbf{Regions of existence for a class of nonlinear diffusion type problems}}
\author{Amit K. Verma$^a$, Mandeep Singh$^b\footnote{Corresponding author Email: $^a$ mandeep04may@yahoo.in}$ , Ravi P. Agarwal$^c$
\\\small{\it{$^a$Department of Mathematics}}\\\small{\it{ Indian Institute of Technology Patna}}\\\small{\it{ Bihta, Patna, 801103, Bihar, India}}
\\\small{\it{$^b$Department of Mathematics}}\\\small{\it{ Jaypee University of Information Technology}}\\\small{\it{Waknaghat, Solan, 173234, Himachal Pradesh, India}}
\\\small{\it{$^c$Department of Mathematics,Texas A\&M, University-Kingsville}}\\\small{\it{700 University Blvd., MSC 172, Kingsville, Texas  78363-8202}}
}
\date{}
\maketitle
\begin{abstract}
The regions of existence are established for a class of two point nonlinear diffusion type boundary value problems (NDBVP)
\begin{eqnarray*}
&&\label{abst-intr-1} -s''(x)-ns'(x)-\frac{m}{x}s'(x)=f(x,s), \qquad m>0,~n\in \mathbb{R},\qquad x\in(0,1),\\
&&\label{abst-intr-2} s'(0)=0, \qquad a_{1}s(1)+a_{2}s'(1)=C,
\end{eqnarray*}
where $a_{1}>0,$ $a_{2}\geq0,~ C\in\mathbb{R}$. These problems arise  very frequently in many branches of engineering, applied mathematics, astronomy, biological system and modern science (see \cite{Gatica1989, GRAY1980, Baxley1991, Chandershekhar1939, Duggan1986, Chambre1952}). By using the concept of upper and lower solutions with monotone constructive technique, we derive some sufficient conditions for existence  in the regions where $\frac{\partial f}{\partial s}\geq0$ and $\frac{\partial f}{\partial s}\leq0$. Theoretical methods are applied for a set of problems which arise in real life.
\end{abstract}

\section{\textbf{Introduction}}
We study the class of two point nonlinear diffusion type boundary value problems (NDBVP),
\begin{eqnarray}
&&\label{eqn-intr-1} -s''(x)-ns'(x)-\frac{m}{x}s'(x)=f(x,s), \qquad m>0,~n\in \mathbb{R},\qquad x\in(0,1),\\
&&\label{eqn-intr-2} s'(0)=0, \qquad a_{1}s(1)+a_{2}s'(1)=C,
\end{eqnarray}
where $a_{1}>0,$ $a_{2}\geq0,~ C\in \mathbb{R}$ and the nonlinear function $f(x,s)$ is continuous and Lipschitz in $D:=\left\{(x, s)\in[0,1]\times \mathbb{R}\right\}$. It is also assumed that $f(x,s)$ is nonsingular and singular with respect to $x\in [0,1]$ and $s$ respectively.

The BVPs \eqref{eqn-intr-1}-\eqref{eqn-intr-2} may be generalised in to the following class of two point nonlinear singular boundary value problems (BVP),
\begin{eqnarray}
&&\label{eqn-intr-3} -\frac{1}{p(x)}\left(p(x)g(x,s)s'\right)'=F(x,s),\\
&&\label{eqn-intr-4} s'(0)=0, \qquad s(1)=B(s(0), s(1), s'(0)).
\end{eqnarray}
Garner and Shivaji \cite{GARNER1990,GARNER1997} studied the diffusion problems defined by \eqref{eqn-intr-3}-\eqref{eqn-intr-4}. 

If we assume $$g(x,s)=e^{\left(nx\right)}, \quad p(x)=x^{m},\quad F(x,s)=e^{\left(nx\right)}f(x,s),$$
and $$B(s(0), s(1), s'(0))= -\frac{{a_{2}}}{a_{1}}s'(1)+\frac{C}{{a_{1}}},$$
then (\ref{eqn-intr-3})-(\ref{eqn-intr-4})  is transformed into the NDBVP (\ref{eqn-intr-1})-(\ref{eqn-intr-2}). Several researchers \cite{Gatica1989, GRAY1980, Baxley1991, Baxley1994, Fink1991, Wang1994, Wang1996, TINEO1994, Ebaid2011, Jia2013} studied similar and other kind of generalized forms of differential equation (\ref{eqn-intr-1}) subject to suitable boundary conditions.

For some special cases, with $m=0,1,2$ and $n=0$, the considered problem  gives the well known Lane-Emden equation, which has been  widely used in astrophysics and modern science such as symmetric solution of shallow membrane caps, heat conduction and source models in the human-head, thermal explosion etc. (see \cite{Chambre1952, Flesch1975, Pandey1996-JDE, Pandey1996-NATMA, Pandey2002, AKV2011, Duggan1986, Dickey1989}). For the numerical solutions of Lane-Emden type problems the reader may refer \cite{akvdt2019, VermaKayenat2018, MSAKVRP2019}. 

We have discussed the region of existence for some real life applications (listed below), which are instances of the given NDBVP (\ref{eqn-intr-1})-(\ref{eqn-intr-2}).

 \begin{itemize}
 \item [$(a)$] The two point Poisson-Boltzmann equation which occurs in the thermal explosion (see \cite{Chambre1952})

\begin{equation}\label{p-1}
\left.
\begin{aligned}
&-s''(x)-\frac{m}{x} s'(x)= \delta e^s,\qquad x\in(0,1),\qquad \delta>0,\\
& s'(0)=0, \qquad s(1)=0
\end{aligned}~~~~~~~~~~\right\}
\end{equation}
 where $m=1~ \text{or}~ 2$ and $n=0$.

 \item [$(b)$] Problem appears for finding the radial stress subject to normal pressure, in a plane circular elastic surface ($m=3$, $n=0$) (see \cite{Dickey1967, TOSAKA1985})
\begin{equation}\label{p-2}
\left.
\begin{aligned}
&-s''(x)-\frac{3}{x} s'(x)=\frac{2}{s^2},\qquad x\in(0,1),\\
&s'(0)=0, \qquad s(1)=\gamma_0, \qquad \gamma_0>0,
\end{aligned}~~~~~~~~~~\right\}
\end{equation}

\item[(c)] The Lane-Emden equation with parameter $\gamma$ (a physical constant), which was derived by Chandrashekhar \cite{Chandershekhar1939} ($m=2$, $n=0$ and $f(x,s)= s^\gamma$). Here we assume $\gamma=5$.
\begin{equation}\label{p-3}
\left.
\begin{aligned}
&-s''(x)-\frac{2}{x}s'(x)=s^5,\qquad x\in(0,1),\\
&\quad s'(0)=0,~~~~s(1)=\sqrt{\frac{3}{4}}.
\end{aligned}\qquad\qquad\right\}
\end{equation}

\item[(d)] The following equation governing the thermal distribution of heat source in the human head (see \cite{Duggan1986})

\begin{equation}\label{p-4}
\left.
\begin{aligned}
&-s''(x)-\frac{2}{x}s'(x)= e^{-s},\qquad x\in(0,1),\\
&\quad s'(0)=0,\quad 2s(1)+s'(1)=0.
\end{aligned}~~~~~~~~~~\right\}
\end{equation}
\end{itemize}
Zhang \cite{YZ1995} observes and justifies that upper and lower solution technique is most promising technique to treat singular boundary value problems. In this work, making use of constructive technique coupled with upper and lower solutions, existence theorems have been established for a form of two point NDBVP in the region $$D:=\left\{(x, s)\in[0,1]\times \mathbb{R} : \alpha_0 \leq s \leq \beta_0\right\}.$$ The functions $\beta_0$ and $\alpha_0$ are called upper and lower solutions and are defined as
\begin{eqnarray}
&&\label{eqn-intr-5-u} -\beta_{0}''(x)-n\beta_{0}'(x)-\frac{m}{x}\beta_{0}'(x)\geq f(x,\beta_{0}), \qquad m>0,~n\in \mathbb{R},\qquad x\in(0,1),\\
&&\label{eqn-intr-6-u} \beta_{0}'(0)=0, \qquad a_{1}\beta_{0}(1)+{a_2}\beta_{0}'(1)\geq C,
\end{eqnarray}
and
\begin{eqnarray}
&&\label{eqn-intr-5-l} -\alpha_{0}''(x)-n\alpha_{0}'(x)-\frac{m}{x}\alpha_{0}'(x)\leq f(x,\alpha_{0}), \qquad m>0,~n\in \mathbb{R},\qquad x\in(0,1),\\
&&\label{eqn-intr-6-l} \alpha_{0}'(0)=0, \qquad a_{1}\alpha_{0}(1)+a_{2}\alpha_{0}'(1)\leq C,
\end{eqnarray}
respectively. The proof of existence results is accomplished by use of following constructive scheme
\begin{eqnarray}
&&\label{eqn-intr-6} -s_{i+1}''(x)-ns_{i+1}'(x)-\frac{m}{x}s_{i+1}'(x)-\lambda s_{i+1}=f(x,s_{i}) - \lambda s_{i},~m>0,~n\in \mathbb{R},~ x\in(0,1)~i=0,1,2,\cdots,\\
&&\label{eqn-intr-7} s_{i+1}'(0)=0, \qquad a_{1}s_{i+1}(1)+a_{2}s_{i+1}'(1)=C,
\end{eqnarray}
where $\lambda$ resembles $\frac{\partial f}{\partial s}$ and it takes both negative and positive values. Under some particular conditions on $f(x,s)$ and $\lambda$, an uniformly convergent (non-increasing monotone) sequence  $\{\beta_{i}\}$, converging to a solution $\widetilde{\beta}(x)$ of (\ref{eqn-intr-1})-(\ref{eqn-intr-2}) in $D$, is generated with the help of initial iterate $\beta_0$. Similarly, making use of the initial iterate $\alpha_{0}$, we also construct the analogous non-decreasing monotone sequence $\{\alpha_{i}\}$ which converges uniformly to a solution ${\widetilde{\alpha}}(x)$ of equation (\ref{eqn-intr-1})-(\ref{eqn-intr-2}) in the region $D$. We conclude that any solution $S(x)$ satisfies ${\widetilde{\alpha}}(x)\leq{S}(x) \leq \widetilde{\beta}(x)$.

This paper is categorized into following sections. In Section \ref{Sec-2}, we construct Green's functions, with the help of Generalized Laguerre polynomial and the hypergeometric function and established Maximum Principle for non-homogeneous linear BVP (\ref{eqn-intr-8})-(\ref{eqn-intr-9}). In Section \ref{Sec-3}, with the support of coupled technique (monotone constructive technique with upper and lower solutions), the regions of existence are established for NDBVP. Finally, the established theory is verified for some real life applications and summarized with concluding remark.
\section{Green's function and its sign} \label{Sec-2}
This section deals with a linear model (non-homogeneous) of the two point NDBVP (\ref{eqn-intr-1})-(\ref{eqn-intr-2}) and the construction of Green's function for the corresponding homogeneous model.

\noindent The iterative technique in this paper gives rise to the following class of linear singular boundary value problem
\begin{eqnarray}
&&\label{eqn-intr-8} -s''(x)-ns'(x)-\frac{m}{x}s'(x)-\lambda s(x)=h(x), \qquad m>0,~n\in \mathbb{R},\qquad x\in(0,1),\\
&&\label{eqn-intr-9} s'(0)=0, \qquad a_{1}s(1)+a_{2}s'(1)=C,
\end{eqnarray}
where $h(x)\in C(I), I=[0,1]$ and $a_{1}>0,$ $a_{2}\geq0,~ C\geq0$.

\noindent Based on $\lambda$, we classify it into the following two subsections:
\subsection{Case-I} For $\lambda=-k^2<0$, where $k\in \mathbb{R}$ (\ref{eqn-intr-8})-(\ref{eqn-intr-9}) reduced to
\begin{eqnarray}
&&\label{eqn-intr-10} -s''(x)-ns'(x)-\frac{m}{x}s'(x)+ k^2 s(x)=0, \qquad m>0,~n\in \mathbb{R},\qquad x\in(0,1),\\
&&\label{eqn-intr-11} s'(0)=0, \qquad a_{1}s(1)+a_{2}s'(1)=0.
\end{eqnarray}
Solution of differential equation (\ref{eqn-intr-10}) is (see \cite{Dinmohammadi2017})
\begin{eqnarray*}
&&\label{eqn-intr-12} s(x)= C_1 e^{-\frac{1}{2}\left(n+\sqrt{4k^2+n^2}\right)x} L_{\alpha}^{m-1}\left(\sqrt{4k^2+n^2}x\right)+C_2 e^{-\frac{1}{2}\left(n+\sqrt{4k^2+n^2}\right)x} U\left({\beta}, m, \sqrt{4k^2+n^2}x\right) ,
\end{eqnarray*}
where $\alpha = \frac{-m\left(n+\sqrt{4k^2 +n^2}\right)}{2\sqrt{4k^2 +n^2}}$, $\beta = \frac{m\left(n+\sqrt{4k^2 +n^2}\right)}{2\sqrt{4k^2 +n^2}}$  and $L_{\alpha}^{m-1}$ and $U(a,b,z)$ are the generalized Laguerre polynomial (related to hydrogen atom wave function) and the hypergeometric function (second linearly independent function) respectively.

Now we define two solutions $u(x)$ and $v(x)$ by
\begin{eqnarray}
&&\label{eqn-intr-13} u(x)= u_m(x),\\
&&\label{eqn-intr-14}v(x)= v_m(x)-A u_m(x),
\end{eqnarray}
where
\begin{eqnarray}
&&u_m(x)= e^{-\frac{1}{2}\left(n+\sqrt{4k^2+n^2}\right)x} L_{\alpha}^{m-1}\left(\sqrt{4k^2+n^2}x\right)\\
&& v_m(x)=e^{-\frac{1}{2}(n+\sqrt{4k^2+n^2})x} U({\beta}, m, \sqrt{4k^2+n^2})x),\\
&& A= \frac{a_{1} v_m(1)+a_{2} v'_{m}(1)}{a_{1} u_m(1)+a_{2} u'_{m}(1)},
\end{eqnarray}
of the differential equation (\ref{eqn-intr-10}) so that $u'(0)=0$ and $a_{1}v(1)+a_{2}v'(1)=0$.
\begin{lemma}\label{Lemma-1}
If $s\in C^2([0, 1])$ is a solution of nonhomogeneous linear two-point BVP (\ref{eqn-intr-8})-(\ref{eqn-intr-9}); then
\begin{eqnarray*}
&& s(x)= \frac{Cu(x)}{a_{1}u(1)+a_{2} u'(1)}-\int_{0}^{1}{G(x,t) h(t)}dt,
\end{eqnarray*}
where the Green's function $G(x,t)$ is defined as
\begin{eqnarray*}
G(x,t)=\frac{1}{W(t)}\left\{
  \begin{array}{ll}
   {u(x)v(t)},~~~x\leq t; \\

   {v(x)u(t)},~~~x\geq t,
  \end{array}
\right.
\end{eqnarray*}
and $W(t)$ is the Wronskian.
\end{lemma}
\textbf{Proof:} The Green's function for the two point BVPs (\ref{eqn-intr-10})-(\ref{eqn-intr-11}) is defined as
\begin{eqnarray*}
G(x,t)=\left\{
  \begin{array}{ll}
   {C_1 u(x)},~~~x\leq t; \\

   {C_2 v(x)},~~~x\geq t.
  \end{array}
\right.
\end{eqnarray*}
Making use of the properties and definition of the Green's function, for $t\in[0, 1]$, we have
\begin{eqnarray*}
&& C_1 u(t) = C_2 v(t)\\
&& C_1 u'(t) - C_2 v'(t)=-1.
\end{eqnarray*}
This gives
\begin{eqnarray*}
&& C_1 = \frac{v(t)}{u(t)v'(t)-v(t)u'(t)} = \frac{v(t)}{W(t)}\\
&& C_2 = \frac{u(t)}{u(t)v'(t)-v(t)u'(t)} = \frac{u(t)}{W(t)}.
\end{eqnarray*}
As $u(x)$ and $v(x)$ satisfy the boundary conditions, which gives
\begin{align*}
W(t) = W(1) e^{n-nt} t^{-m},
\end{align*}
where $W(1)$ depends on $n$, $m$ and $k$.
Finally, we get the Green's function $G(x,t)$ for linear two point BVP (\ref{eqn-intr-10})-(\ref{eqn-intr-11}).

The linear BVP (\ref{eqn-intr-10})-(\ref{eqn-intr-11}) is equivalent to the integral equation
\begin{eqnarray*}
&& s(x)= \bar{s}-\int_{0}^{1}{G(x,t) h(t)}dt,
\end{eqnarray*}
where $\bar{s}$ is the solution of nonhomogeneous linear BVP (\ref{eqn-intr-8})-(\ref{eqn-intr-9})
\begin{eqnarray*}
&& \bar{s}= C_1 u(x)+C_2v(x).
\end{eqnarray*}
Now from the boundary conditions $u'(0)=0$ and $a_{1}v(1)+a_{2}v'(1)=C$, we get
\begin{eqnarray*}
&& \bar{s}= \frac{Cu(x)}{a_{1}u(1)+a_{2} u'(1)}.
\end{eqnarray*}
Hence, the linear BVP (\ref{eqn-intr-8})-(\ref{eqn-intr-9}) is equivalent to
\begin{eqnarray*}
&& s(x)= \frac{Cu(x)}{a_{1}u(1)+a_{2} u'(1)}-\int_{0}^{1}{G(x,t) h(t)}dt.
\end{eqnarray*}
\begin{lemma}\label{Lemma-2}
The sign of the Green's function
\begin{eqnarray*}
G(x,t)=\frac{1}{W(t)}\left\{
  \begin{array}{ll}
   {u(x)v(t)},~~~x\leq t; \\

   {v(x)u(t)},~~~x\geq t,
  \end{array}
\right.
\end{eqnarray*}
 for linear singular boundary value problem (\ref{eqn-intr-10})-(\ref{eqn-intr-11}) is non-positive.
\end{lemma}
\textbf{Proof.} The solutions $u(x)$ and $v(x)$, for $a_{1}>0$, $a_{2}\geq 0$, $m>0$ and $n\in R$ (defined by (\ref{eqn-intr-13})-(\ref{eqn-intr-14})) of differential equation (\ref{eqn-intr-10}) satisfy the following properties.
\begin{itemize}
\item The solution $v(x)$ is a positive (i.e., $v(x)>0$), non-increasing function for $x\in[0,1]$ and unbounded at the origin.
\item If $W(1)$ is a positive for given $m$, $n$ and $k$, then $u(x)$ is a negative (i.e., $u(x)<0$), non-increasing function of $x\in [0, 1]$.
\item Similarly, if $W(1)$ is negative for given $m$, $n$ and $k$, then $u(x)$ is a positive (i.e., $u(x)>0$), non-decreasing function of $x\in [0, 1]$.
\end{itemize}
Using of the above facts, we can easily conclude that $G(x,t)\leq0.$
\subsection{Case II} Similar to Case I, for $\lambda=k^2$ we consider the following homogeneous linear BVP
\begin{eqnarray}
&&\label{eqn-intr-15} -s''(x)-ns'(x)-\frac{m}{x}s'(x)-k^2 s(x)=0, \qquad m>0,~n\in \mathbb{R},\qquad x\in(0,1),\\
&&\label{eqn-intr-16} s'(0)=0, \qquad a_{1}s(1)+a_{2}s'(1)=0.
\end{eqnarray}
The two solutions $(\bar{u}(x))$ and $(\bar{v}(x))$ of differential equation (\ref{eqn-intr-15}) are
\begin{eqnarray}
&&\label{eqn-intr-17} \bar{u}(x)= \bar{u}_m(x),\\
&&\label{eqn-intr-18}\bar{v}(x)= \bar{v}_m(x)-A \bar{u}_m(x),
\end{eqnarray}
where
\begin{eqnarray}
&&\bar{u}_m(x)= e^{-\frac{1}{2}\left(n+\sqrt{n^2-4k^2}\right)x} L_{\bar{\alpha}}^{m-1}\left(\sqrt{n^2-4k^2}x\right),\\
&& \bar{v}_m(x)=e^{-\frac{1}{2}\left(n+\sqrt{n^2-4k^2}\right)x} U\left({\beta}, m, \sqrt{n^2-4k^2}x\right),\\
&& S= \frac{a_{1} \bar{v}_m(1)+a_{2} \bar{v}'_{m}(1)}{a_{1} \bar{u}_m(1)+a_{2} \bar{u}'_{m}(1)},\\
&& \bar{\alpha} = \frac{-m(n+\sqrt{n^2-4k^2})}{2\sqrt{n^2-4k^2}},\\
&& \bar{\beta} = \frac{m(n+\sqrt{n^2-4k^2})}{2\sqrt{n^2-4k^2}},
\end{eqnarray}
and so $\bar{u}'(0)=0$ and $a_{1} \bar{v}(1)+a_{2}\bar{v}'(1)=0$.

\begin{lemma}\label{Lemma-3}
If $s\in C^2([0,1])$ is a solution of nonhomogeneous linear two point boundary value problem (\ref{eqn-intr-8})-(\ref{eqn-intr-9}), then
\begin{eqnarray*}
&& s(x)= \frac{C\bar{u}(x)}{a_{1}\bar{u}(1)+a_{2} \bar{u}'(1)}-\int_{0}^{1}{G(x,t) h(t)}dt,
\end{eqnarray*}
where the Green's function $G(x,t)$ is defined as
\begin{eqnarray*}
G(x,t)=\frac{1}{W(t)}\left\{
  \begin{array}{ll}
   {\bar{u}(x)\bar{v}(t)},~~~x\leq t, \\

   {\bar{v}(x)\bar{u}(t)},~~~x\geq t,
  \end{array}
\right.
\end{eqnarray*}
and $W(t)$ is the Wronskian.
\end{lemma}
\textbf{Proof:} Proof is similar as Lemma \ref{Lemma-1}.

\begin{lemma}\label{Lemma-2g} Let $\lambda=k^2<\lambda_0$, where $\lambda_0$ is the first positive eigenvalue of \eqref{eqn-intr-10}-\eqref{eqn-intr-11}. The sign of the Green's function $G(x,t)$ defined in Lemma \ref{Lemma-3} for linear singular boundary value problems (\ref{eqn-intr-10})-(\ref{eqn-intr-11}) is non-positive.
\end{lemma}
\textbf{Proof:} Proof is similar as Lemma \ref{Lemma-2}.

\subsection{Maximum Principle}
\begin{proposition}\label{Proposition-1}
Let $s(x)$ be the solution of linear BVP (\ref{eqn-intr-8})-(\ref{eqn-intr-9}) and let $\lambda<\lambda_0$. If $b$ is non-negative and $h(x)$ is continuous in $[0,1]$ such that $h(x)\geq 0$, then $s(x)\geq 0$.
\end{proposition}
\textbf{Proof:} Proof is similar as Proposition $3.1$ of \cite{Verma2015}.

\section{Nonlinear Diffusion BVP : The Region of Existence} \label{Sec-3}
In this section, making use of Maximum Principle we generate monotone sequences with the help of coupled technique, which, finally establishes the existence results (under certain conditions) for the two point NDBVP.
\begin{eqnarray*}
&&-s''(x)-ns'(x)-\frac{m}{x}s'(x)=f(x,s), \qquad m>0,~n\in \mathbb{R},\qquad x\in(0,1),\\
&& s'(0)=0, \qquad a_{1}s(1)+a_{2}s'(1)=C,
\end{eqnarray*}
where $a_{1}>0,$ $a_{2}\geq0$ and $C\in\mathbb{R}$.

Now we define one sided Lipschitz condition which was introduced by L. F. Shampine et al. \cite{SHAMPINE1970}.
\begin{definition}
A function $f(x,y)$ is said to satisfy one sided Lipschitz condition if $y\leq w$ implies $f(x,y)-f(x,w)\leq L (w-y)$. Here $L$ may be referred as one sided Lipschitz constant.
\end{definition}

The following theorem is one of two important results of this paper. 
\begin{theorem}\label{theorem-1}
The two point NDBVP (\ref{eqn-intr-1})-(\ref{eqn-intr-2}) has at the minimum one real solution in the region $$D = \{(x,s) \in [0,1] \times \mathbb{R} :\alpha_0 \leq s \leq \beta_0\}$$ under the following assumptions
\begin{itemize}
\item [$(a)$] there exist upper ($\beta_0$) and lower ($\alpha_0$) solutions in $C^2[0,1]$, defined by (\ref{eqn-intr-5-u})-(\ref{eqn-intr-6-u}) and (\ref{eqn-intr-5-l})-(\ref{eqn-intr-6-l}) respectively, such that $\beta_0\geq \alpha_0$;
\item [$(b)$] the source function $f:D\rightarrow\mathbb{R}$ is continuous on $D$;
\item [$(c)$] $\exists$ one-sided Lipschitz constant $ L \geq 0$ such that $\forall$ $(x,w_{1}),(x,w_{2})\in D$
\begin{eqnarray*}
w_{2}\leq w_{1} \Longrightarrow  f(x,w_{1})- f(x,w_{2}) \geq - L(w_{1} - w_{2}).
\end{eqnarray*}
\end{itemize}
Further if there is a $\lambda$ (constant) such that $\lambda\leq0$, and  $ L+\lambda \leq 0$ then the monotonically (non-increasing) sequences $\{\beta_i\}$ with initial approximate $\beta_0$, generated by
\begin{eqnarray}
&&\label{eqn-MIT-1} -\beta_{i+1}''(x)-n\beta_{i+1}'(x)-\frac{m}{x}\beta_{i+1}'(x)-\lambda \beta_{i+1}=F(x,\beta_{i}), \qquad m>0,~n\in \mathbb{R},\qquad x\in(0,1),\\
&&\label{eqn-MIT-2} \beta_{i+1}'(0)=0, \qquad a_{1}\beta_{i+1}(1)+a_{2}\beta_{i+1}'(1)=C,
\end{eqnarray}
where $F(x,\beta_{i})=f(x,\beta_{i})-\lambda \beta_{i}$, converges uniformly to a solution ${{\widetilde{\beta}}(x)}$ of (\ref{eqn-intr-1})-(\ref{eqn-intr-2}). Similarly using initial approximate $\alpha_0$, leads to a monotonically (non-decreasing) sequences $\{\alpha_i\}$ converges uniformly towards a solution ${\widetilde{\alpha}}(x)$. Any solution of the considered problem  (\ref{eqn-intr-1})-(\ref{eqn-intr-2}), (say ${S}(x)\in D$) must satisfy
\begin{eqnarray*}
{\widetilde{\alpha}}(x)\leq{S}(x) \leq \widetilde{\beta}(x).
\end{eqnarray*}
\end{theorem}
\textbf{Proof.} Making use of (\ref{eqn-intr-5-u})-(\ref{eqn-intr-6-u}) and (\ref{eqn-MIT-1})-(\ref{eqn-MIT-2}) (for $i=0$), we get
\begin{eqnarray*}
&& -(\beta_0-\beta_1)''(x)-(n+\frac{m}{x})(\beta_0-\beta_1)'(x)-\lambda (\beta_0-\beta_1)(x)\geq 0,\\
&&~(\beta_0-\beta_1)'(0)=0, ~~~~a_{1}(\beta_0-\beta_1)(1)+a_{2}(\beta_0-\beta_1)'(1)\geq 0.
\end{eqnarray*}
As $h(x)\geq0$ and $b\geq 0$, by using the Maximum principle, Proposition \ref{Proposition-1}, we have $\beta_0 \geq \beta_1$.

By using the concept of one sided Lipschitz condition and the assumption $L+\lambda \leq0 $, from equation (\ref{eqn-MIT-1})-(\ref{eqn-MIT-2}) we get
\begin{eqnarray*}
&&-\beta_{i+1}''(x)-(n+\frac{m}{x})\beta_{i+1}'(x)-\lambda \beta_{i+1} \geq (L+\lambda)(\beta_{i+1}-\beta_{i})+f(x,\beta_{i+1}).
\end{eqnarray*}
And if $(\beta_{i}\geq \beta_{i+1})$, then
\begin{eqnarray}
\label{Th-1-eqn-1}&&-\beta_{i+1}''(x)-(n+\frac{m}{x})\beta_{i+1}'(x)-\lambda \beta_{i+1}(x) \geq f(t,\beta_{i+1}).
\end{eqnarray}
Since $\beta_0 \geq \beta_1$, then from equation (\ref{Th-1-eqn-1}) (for $n=0$) and (\ref{eqn-MIT-1})-(\ref{eqn-MIT-2}) (for $i=1$) we get
\begin{eqnarray*}
&&-(\beta_1-\beta_2)''(x)-(n+\frac{m}{x})(\beta_1-\beta_2)'(x)-\lambda (\beta_1-\beta_2)(x)\geq 0,\\
&&~(\beta_1-\beta_2)'(0)=0, ~~~~a_{1}(\beta_1-\beta_2)(1)+a_{2}(\beta_1-\beta_2)'(1)\geq 0.
\end{eqnarray*}
From the Maximum principle, Proposition \ref{Proposition-1}, we have $\beta_1\geq \beta_2$.

\noindent Now with the help of lower solution (\ref{eqn-intr-5-l})-(\ref{eqn-intr-6-l}) and (\ref{eqn-MIT-1})-(\ref{eqn-MIT-2}) (for $i=0$)
\begin{eqnarray*}
&&-(\beta_1-\alpha_0)''(x)-(n+\frac{m}{x})(\beta_1-\alpha_0)'(x)-\lambda (\beta_1-\alpha_0)(x)\geq 0,\\
&&(\beta_1-\alpha_0)'(0)=0\quad a_{1}(\beta_1-\alpha_0)(1)+a_{2}(\beta_1-\alpha_0)'(1)\geq 0.
\end{eqnarray*}
Thus, $\beta_1\geq \alpha_0$ follows from the Maximum Principle.

\noindent To generate monotonic sequences, we assume $\beta_{i}\geq \beta_{i+1}$, $\beta_{i+1}\geq \alpha_{0}$ and show that $\beta_{i+1}\geq \beta_{i+2}$,  $\beta_{i+2}\geq \alpha_{0}$ for all $i$.

From equations (\ref{eqn-MIT-1})-(\ref{eqn-MIT-2}) (for $(i+1)^{th}$) and (\ref{Th-1-eqn-1}), we get
\begin{eqnarray*}
&&-(\beta_{i+1}-\beta_{i+2})''(x)-(n+\frac{m}{x})(\beta_{i+1}-\beta_{i+2})'(x)-\lambda (\beta_{i+1}-\beta_{i+2})(x)\geq 0,\\
&&~~(\beta_{i+1}-\beta_{i+2})'(0)=0,~~~~a_{1}(\beta_{i+1}-\beta_{i+2})(1)+a_{2}(\beta_{i+1}-\beta_{i+2})'(1)\geq 0,
\end{eqnarray*}
and hence from the Maximum Principle, we have $\beta_{i+1}\geq \beta_{i+2}$.

\noindent From equations (\ref{eqn-MIT-1})-(\ref{eqn-MIT-2}) (for $(i+1)^{th}$) and (\ref{Th-1-eqn-1}) we get,
\begin{eqnarray*}
&&-(\beta_{i+2}-\alpha_{0})''(x)-(n+\frac{m}{x})(\beta_{i+2}-\alpha_{0})'(x)-\lambda (\beta_{i+2}-\alpha_{0})(x)\geq 0,\\
&&(\beta_{i+2}-\alpha_{0})'(0)=0,~~~~a_{1}(\beta_{i+2}-\alpha_{0})(1)+a_{2}(\beta_{i+2}-\alpha_{0})'(1)\geq0.
\end{eqnarray*}
Then, from the Maximum Principle, $\beta_{i+2}\geq \alpha_{0}$ and hence we have
\begin{eqnarray*}
\alpha_{0}\leq \ldots \leq \beta_{i+1}\leq \beta_{i}\leq \ldots \leq \beta_{2}\leq \beta_{1}\leq \beta_{0}.
\end{eqnarray*}
Starting with $\alpha_0$ it is easy to get
\begin{eqnarray*}
\alpha_{0} \leq \alpha_{1}\leq \alpha_{2}\leq  \ldots \leq \alpha_{i}\leq \alpha_{i+1}\leq  \ldots \leq \beta_{0}.
\end{eqnarray*}
Finally we show that $\beta_{i}\geq \alpha_{i}$ for all $i$. For this, by assuming $\beta_{i}\geq \alpha_{i}$, we show that $\beta_{i+1}\geq \alpha_{i+1}$. From equation (\ref{eqn-MIT-1})-(\ref{eqn-MIT-2}), we get
\begin{eqnarray*}
&&-(\beta_{i+1}-\alpha_{i+1})''(x)-(n+\frac{m}{x})(\beta_{i+1}-\alpha_{i+1})'(x)-\lambda (\beta_{i+1}-\alpha_{i+1})(x)\geq 0,\\
&&~~(\beta_{i+1}-\alpha_{i+1})'(0)=0,~~~~a_{1}(\beta_{i+1}-\alpha_{i+1})(1)+a_{2}(\beta_{i+1}-\alpha_{i+1})'(1)\geq 0.
\end{eqnarray*}
Hence from the Maximum Principle, $\beta_{i+1}\leq \alpha_{i+1}$. Thus, we have
\begin{eqnarray*}
\alpha_{0}\leq \alpha_{1}\leq \alpha_{2}\leq  \ldots \leq \alpha_{i}\leq \alpha_{i+1}\leq  \ldots \leq \beta_{i+1}\leq \beta_{i}\leq \ldots \leq \beta_{2}\leq \beta_{1}\leq \beta_{0}.
\end{eqnarray*}
This gives bounded non-decreasing $\{(\alpha_{i})\}$ and non-increasing  $\{(\beta_{i})\}$ sequences. Now by applying Dini's theorem, we conclude that the governing sequences converge uniformly. Let $\widetilde{\beta}(x)=\displaystyle\lim_{i\to\infty}\beta_{i}(x) $ and $\widetilde{\alpha}(x)=\displaystyle\lim_{i\to\infty}\alpha_{i}(x)$.

By Lemma \ref{Lemma-1}, the solution $\beta_{i+1}$ of equation (\ref{eqn-MIT-1})-(\ref{eqn-MIT-2}) is given by
\begin{eqnarray*}
&& \beta_{i+1}= \frac{Cu(x)}{a_{1}u(1)+a_{2} u'(1)}-\int_{0}^{1}{G(x,t)(f(t,\beta_{i})-\lambda \beta{i})}dt.
\end{eqnarray*}
Then, making use of the Lebesgue's dominated convergence theorem, as $i \longrightarrow \infty$, we get
\begin{eqnarray*}
&&\widetilde{\beta}(x)= \frac{Cu(x)}{a_{1}u(1)+a_{2} u'(1)}- {\int_0}^1{ G(x,t)  (f(t,\widetilde{\beta})-\lambda \widetilde{\beta})dt}.
\end{eqnarray*}
This is the solution of NDBVP (\ref{eqn-intr-1})-(\ref{eqn-intr-2}). Similar equation can be arise for lower solution also. Any result $S(x)\in D$ of (\ref{eqn-MIT-1})-(\ref{eqn-MIT-2}) can take the part of $\beta_0(x)$, hence $S(x)\geq \widetilde{\alpha}(x)$ and in the similar way, we conclude that $S(x)\leq \widetilde{\beta}(x)$.

Now we give another important results in the form of following theorem.
\begin{theorem}\label{theorem-2}
The two point NDBVP (\ref{eqn-intr-1})-(\ref{eqn-intr-2}) has at the minimum one real solution in the region $$\widetilde{D}:=\{(x,s) \in [0,1] \times \mathbb{R} :\alpha_0 \leq s \leq \beta_0 \}$$ under the following assumptions
\begin{itemize}
\item [$(a)$] there exist upper ($\beta_0$) and lower ($\alpha_0$) solutions in $C^2[0,1]$, defined by (\ref{eqn-intr-5-u})-(\ref{eqn-intr-6-u}) and (\ref{eqn-intr-5-l})-(\ref{eqn-intr-6-l}) respectively, such that $\beta_0\geq \alpha_0$;
\item [$(b)$] the source function $f: \widetilde{D} \rightarrow \mathbb{R}$ is continuous on $\widetilde{D}$;
\item [$(c)$] $\exists$ one-sided Lipschitz constant $ L\geq 0$ such that $\forall$ $(x,w_{1}),(x,w_{2})\in \widetilde{D}$
\begin{eqnarray*}
w_{2}\leq w_{1} \Longrightarrow  f(x,w_{1})- f(x,w_{2}) \geq L(w_{1} - w_{2}).
\end{eqnarray*}
\end{itemize}
Further if there is a  $\lambda$ (constant) such that $0<\lambda<\lambda_0$ where $\lambda_0$ is the first eigenvalue of \eqref{eqn-intr-10}-\eqref{eqn-intr-11} and $ \lambda-L \leq 0$ then the monotonically (non-increasing) sequences $\{\beta_i\}$, with initial approximate $\beta_0$ generated by
\begin{eqnarray*}
&& -\beta_{i+1}''(x)-n\beta_{i+1}'(x)-\frac{m}{x}\beta_{i+1}'(x)-\lambda \beta_{i+1}=F(x,\beta_{i}), \qquad m>0,~n\in \mathbb{R},\qquad x\in(0,1),\\
&&\beta_{i+1}'(0)=0, \qquad a_1\beta_{i+1}(1)+a_2\beta_{i+1}'(1)=C,
\end{eqnarray*}
where $F(x,\beta_{i})=f(x,\beta_{i})-\lambda \beta_{i}$,  converges uniformly to a solution ${{\widetilde{\beta}}(x)}$ of (\ref{eqn-intr-1})-(\ref{eqn-intr-2}). Similarly, using initial approximate $\alpha_0$, leads to a monotonically (non-decreasing) sequences $\{\alpha_i\}$ converges uniformly towards a solution ${\widetilde{\alpha}}(x)$.  Any solution of  the considered problem  (\ref{eqn-intr-1})-(\ref{eqn-intr-2}), (say ${S}(x)\in D$) must satisfy
\begin{eqnarray*}
{\widetilde{\alpha}}(x)\leq{S}(x) \leq \widetilde{\beta}(x).
\end{eqnarray*}
\end{theorem}
\textbf{Proof.}  Proof follows similar analysis to the proof of Theorem \ref{theorem-1}.

\section{Real Life Applications}
This section deals with some two point NDBVP that arise in various real life problems occurring in sciences and engineering. We verify our results and show that under some sufficient conditions, it is possible that we can generate monotone sequences that converges to solution of NDBVP.
\begin{example} Consider the  two point nonlinear diffusion BVP
\begin{eqnarray}
&&-s''(x)-s'(x)-\frac{1}{x} s'(x)=1- 2e^{s}\\
&& s'(0)=0, \qquad s(1)+s'(1)=1.
\end{eqnarray}
\end{example}
Here $f(x, s)=1-2e^{s}$. $a_{1}=1,~a_{2}=1$ and $C=1$. This problem has lower solution $\alpha_0=-1$ and upper solution $\beta_0=1$. The nonlinear source function $f(x, s)=1-2e^{s}$ is Lipschitz and continuous for all $s$ in domain $D$. The one sided Lipschitz constant is $L \geq2e$. Now by using Theorem \ref{theorem-1}, we can conclude that for some $\lambda=-k^2<0$, such that $\lambda\leq-2e$, we can generate two monotonically (non-increasing and non-decreasing) sequences, which converge to the solutions of given nonlinear singular BVP. The region of existence is given by $$\{(x,s):0\leq x\leq 1, -1\leq s\leq 1\}.$$

\begin{example} Consider the following nonlinear diffusion BVP
\begin{eqnarray}
&&-s''(x)-10s'(x)-\frac{5}{x} s'(x)=\frac{3}{4}{s}\\
&& s'(0)=0, \qquad 6s(1)+ s'(1)=1.
\end{eqnarray}
\end{example}
Here $f(x, y)=\frac{3}{4}{s(x)}$. $a_{1}=6,~a_{2}=1$ and $C=1$. This problem has  $\alpha_0=0$ and $\beta_0=\frac{2-x^2}{3}$ as lower and upper solutions respectively. The nonlinear source function $f(x, s)=\frac{3}{4}{s(x)}$ is Lipschitz and continuous for all $s$ in domain $D$. The one sided Lipschitz constant is $L\leq\frac{3}{4}$. Now by using Theorem \ref{theorem-2}, we can conclude that for some $\lambda=k^2>0$, such that $\lambda\leq L$, we can generate two monotonically (non-increasing and non-decreasing) sequences, which converge to the solutions of given nonlinear singular BVP. The region of existence is given by $$\left\{(x,s):0\leq x\leq 1, 0\leq s\leq \frac{2-x^2}{3}\right\}.$$

\begin{example} Thermal explosion in cylindrical vessel (from equation (\ref{p-1}) with $m=1$ and $\delta=\frac{1}{4}$) 
\begin{eqnarray}
&&-s''(x)-\frac{1}{x} s'(x)=\frac{1}{4} e^s\\
&& s'(0)=0, \qquad s(1)=0.
\end{eqnarray}
\end{example}
Here $f(x, s)=\frac{1}{4} e^s$, $a_{1}=1,~a_{2}=0$ and $C=1$. This problem has  $\alpha_0=0$ and $\beta_0=\frac{3-x^2}{4}$ as lower and upper solutions respectively. The nonlinear source function $f(x, s)=\frac{1}{4} e^s$ is Lipschitz and continuous for all $s$ in domain $D$. The one sided Lipschitz constant is $L\leq\frac{1}{4}$. Now by using Theorem \ref{theorem-2}, we can conclude that for some $\lambda=k^2>0$, such that $\lambda\leq L$, we can generate two monotonically (non-increasing and non-decreasing) sequences, which converge to the solutions of given nonlinear singular BVP. The region of existence is given by $$\left\{(x,s):0\leq x\leq 1, 0\leq s\leq \frac{3-x^2}{4}\right\}.$$

\begin{example} The radial stress in a plane circular elastic surface subject to normal pressure (see equation (\ref{p-2}))
\end{example}
Here $f(x, s)=\frac{2}{s^2}$, $a_{1}=1,~a_{2}=0$, $C=1$ and $\gamma_0=2$. This problem has  $\alpha_0=2$ and $\beta_0=2+\frac{1-x^2}{9}$ as lower and upper solutions respectively. The nonlinear source function $f(x, s)$  is Lipschitz in $s$ and also continuous for all $s$ in domain $D:=\{(x,s) \in [0,1] \times R :\alpha_0 \leq s \leq \beta_0\}$. The one sided Lipschitz constant is $L\geq\frac{19}{36}$. Now by using Theorem \ref{theorem-1}, we can conclude that for some $\lambda=-k^2<0$, such that $\lambda\leq -L$, we can generate two monotonically (non-increasing and non-decreasing) sequences, which converge to the solutions of given nonlinear singular BVP. The region of existence is given by $$\left\{(x,s):0\leq x\leq 1, 2\leq s\leq 2+\frac{1-x^2}{9}\right\}.$$

\begin{example} Lane-Emden equation (see equation (\ref{p-3}))

\end{example}
Here $f(x, s)=s^5$, $a_{1}=1,~a_{2}=0$ and $C=\sqrt{\frac{3}{4}}$. This problem has  $\alpha_0=\frac{3}{4}$ and $\beta_0=\frac{1}{\sqrt{0.7+\frac{x^2}{2}}}$ as lower and upper solutions respectively. The nonlinear source function $f(x, s)=s^5$ is Lipschitz and continuous for all $y$ in domain $D$. The one sided Lipschitz constant is $L\leq1.58203125$. Now by using Theorem \ref{theorem-2}, we can conclude that for some $\lambda=k^2>0$, such that $\lambda\leq L$, we can generate two monotonically (non-increasing and non-decreasing) sequences, which converge to the solutions of given nonlinear singular BVP. The region of existence is given by $$\left\{(x,s):0\leq x\leq 1, \frac{3}{4}\leq s\leq \frac{1}{\sqrt{0.7+\frac{x^2}{2}}}\right\}.$$

\begin{example} Thermal distribution in the human heads (see equation (\ref{p-4})).

\end{example}
Here $f(x, s)=e^{-s}$, $a_{1}=2,~a_{2}=1$ and $C=0$. This problem has  $\alpha_0=0$ and $\beta_0={2-x^2}$ as lower and upper solutions respectively. The nonlinear source function $f(x, y)=e^{-s}$ is Lipschitz and continuous for all $s$ in domain $D$. The one sided Lipschitz constant is $L\geq e^{2}$. Now by using Theorem \ref{theorem-1}, we can conclude that for some $\lambda=-k^2<0$, such that $\lambda\leq - L$, we can generate two monotonically (non-increasing and non-decreasing) sequences, which converge to the solutions of given nonlinear singular BVP. The region of existence is given by $$\{(x,s):0\leq x\leq 1, 0\leq s\leq 2-x^2\}.$$

\section{Conclusion} In this work, a coupled technique (monotone constructive technique with upper and lower solutions) has been used to established the existence results for a class of nonlinear diffusion type boundary value problems. Green's functions are constructed with the help of Generalized Laguerre polynomial and the hypergeometric function. The established techniques are also tested on some real life problems which are modeled by the given nonlinear diffusion BVP such as thermal explosion in cylindrical vessel, thermal distribution in the human heads and astrophysics. Regions for existence of solutions are computed. Since choice of Initial guesses are not unique, the region of existence computed here is not unique and it can be improved provided choice of initial guesses are better. 

\section{Acknowledgement} We are thankful to the anonymous reviewers for their valuable time and suggestion which helped authors a lot and improved the papers in all respect.
\bibliography{ref}
\bibliographystyle{plain}

\end{document}